\documentclass[reqno,12pt]{amsart}
\usepackage{amsmath,amssymb,latexsym}	
\numberwithin{table}{section}

\newcommand{\Sp}{\mathrm{Sp}}

\newcommand{\GL}{\mathrm{GL}}
\newcommand{\SL}{\mathrm{SL}}

\newcommand{\R}{\mathbb{R}}
\newcommand{\Q}{\mathbb{Q}}
\newcommand{\Z}{\mathbb{Z}}

\newcommand{\C}{\mathbb{C}}
\newcommand{\e}{\epsilon}

\newtheorem{theorem}{Theorem}[section]

\theoremstyle{definition}

\newtheorem{rmk}[theorem]{Remark}

\author{Sandip  Singh}
\address{Institut f\"ur Mathematik,
Johannes Gutenberg-Universit\"at Mainz
55099 Mainz/Germany}
\email{singhs@uni-mainz.de}  
\subjclass[2010]{Primary: 22E40;  Secondary: 32S40;  33C80}  \keywords{Monodromy representation, Hypergeometric group, Symplectic group, Calabi-Yau threefolds}

\begin{document}     \title[Monodromy Groups associated to Calabi-Yau threefolds]{Arithmeticity of Four Hypergeometric Monodromy Groups associated to Calabi-Yau threefolds}

\vskip 5mm
\begin{abstract} In \cite{SV}, we show that 3 of the 14 hypergeometric monodromy groups associated to Calabi-Yau threefolds, are arithmetic. Brav-Thomas (in \cite{BT}) show that 7 of the remaining 11, are thin. In this article, we settle the arithmeticity problem for the 14 monodromy groups, by showing that, the remaining 4 are arithmetic.
\end{abstract}
\maketitle

\section{Introduction}
For $\theta=z\frac{d}{dz}$, we write the differential operator
\begin{align*}
D&=D(\alpha;\beta)=D(\alpha_1,\ldots,\alpha_n;\beta_1,\ldots,\beta_n)\nonumber\\
&=(\theta+\beta_1-1)\cdots(\theta+\beta_n-1)-z(\theta+\alpha_1)\cdots(\theta+\alpha_n)
\end{align*}
for $\alpha_1,\ldots,\alpha_n,\beta_1,\ldots,\beta_n\in\C$, and consider the {\it hypergeometric} differential equation
\begin{eqnarray}\label{introdifferentialequation}
D(\alpha;\beta)w=0
\end{eqnarray}
on $\mathbb{P}^1(\C)$ with {\it regular singularities} at the points $0, 1$ and $\infty$, and regular elsewhere.

The fundamental group of $\mathbb{P}^1(\C)\backslash\{0,1,\infty\}$ acts on the solution space of the differential equation (\ref{introdifferentialequation}), and we get a representation $M(\alpha;\beta)$ of $\pi_1\left(\mathbb{P}^1(\C)\backslash\{0,1,\infty\}\right)$ inside $\GL_n(\C)$, called {\it monodromy}; and the {\it monodromy group} of the hypergeometric equation (\ref{introdifferentialequation}) is the image of this map, i.e. the subgroup of $\GL_n(\C)$ generated by the monodromy matrices $M(\alpha;\beta)(h_0)$, $M(\alpha;\beta)(h_1)$, $M(\alpha;\beta)(h_\infty)$, where $h_0, h_1, h_\infty$ (loops around $0, 1, \infty$ resp.) are the generators of $\pi_1\left(\mathbb{P}^1(\C)\backslash\{0,1,\infty\}\right)$ with a single relation $h_\infty h_1 h_0=1$.

By a theorem of Levelt (\cite{Le}; cf. \cite[Theorem 3.5]{BH}), if $\alpha_1,\alpha_2,\ldots,\alpha_n$, $\beta_1,\beta_2,\ldots,$ $\beta_n\in\C$ such that $\alpha_j-\beta_k\not\in\Z$, for all $j,k=1,2,\ldots,n$,  then the monodromy group of the hypergeometric differential equation (\ref{introdifferentialequation}) is  (up to conjugation in $\GL_n(\C)$)  a subgroup of $\GL_n(\C)$ generated by the {\it companion matrices} $A$ and $B$ of  $$f(X)=\prod_{j=1}^{n}(X-{\rm{e}^{2\pi i\alpha_j}})\quad\mbox{ and }\quad g(X)=\prod_{j=1}^{n}(X-{\rm{e}^{2\pi i\beta_j}})$$ resp., and the monodromy is defined by  $h_\infty\mapsto A$, $h_0\mapsto B^{-1}$, $h_1\mapsto A^{-1}B$.

Therefore, to study the monodromy groups of $n$-order hypergeometric differential equations, it is enough to study the subgroups of $\GL_n(\C)$ generated by the companion matrices of pairs of degree $n$ polynomials in $\C[X]$, which do not have any common root in $\C$. Here we concentrate only on the case $n=4$ and $\alpha,\beta\in\Q^4$.

Let $f,g\in\Z[X]$ be a pair of degree four polynomials which are product of cyclotomic polynomials, do not have any common root in $\C$, $f(0)=g(0)=1$ and form a primitive pair i.e., $f(X)\neq f_1(X^k)$ and $g(X)\neq g_1(X^k)$, for any $k\geq2$ and $f_1,g_1\in\Z[X]$. Now, form the companion matrices $A, B$ of $f, g$ resp., and consider the subgroup $\Gamma\subset\SL_4(\Z)$ generated by $A$ and $B$. It follows from \cite{BH} that $\Gamma$ preserves a non-degenerate integral symplectic form $\Omega$ on $\Z^4$ and $\Gamma\subset\Sp_4(\Omega)(\Z)$ is Zariski dense. The case when $\Gamma\subset\Sp_4(\Omega)(\Z)$ has finite index, it is called {\it arithmetic}, and {\it thin} in other case \cite{S}.

We now take particular examples where $f=(X-1)^4$ (i.e. the local monodromy is maximally unipotent at $\infty$ and $\alpha=(0,0,0,0)$). Then, it turns out (for a reference, see \cite{AvEvSZ}, \cite{YYCE} and \cite{DoMo}) that, this monodromy is same as the monodromy of $\pi_1\left(\mathbb{P}^1(\C)\backslash\{0,1,\infty\}\right)$ on certain pieces of $\mathrm{H}^3$ of the fibre of a family $\{Y_t : t\in\mathbb{P}^1(\C)\backslash\{0,1,\infty\}\}$ of Calabi-Yau threefolds, provided $f,g$ satisfy the conditions that, $f(X)=(X-1)^4$ and $g(X)$ is the product of cyclotomic polynomials such that $g(1)\neq 0$, $f(0)=g(0)=1$, and $f,g$ form a primitive pair. There are precisely $14$ such examples, which have been listed in \cite{AvEvSZ}, \cite{YYCE}, \cite{DoMo}, \cite{SV} and Table \ref{table:calabiyau} of this article.  It is then of interest to know whether the associated monodromy group $\Gamma$ is arithmetic. 

In \cite{SV}, we prove the {\it arithmeticity} of 3 monodromy groups associated to Examples \ref{arithmeticYYCE-1}, \ref{arithmeticYYCE-2}, \ref{arithmeticYYCE-3} of  Table \ref{table:calabiyau}. Brav-Thomas (in \cite{BT}) prove the {\it thinness} of 7 monodromy groups associated to Examples \ref{brav1}, \ref{brav2}, \ref{brav3}, \ref{brav4}, \ref{brav5}, \ref{brav6}, \ref{brav7} of Table \ref{table:calabiyau}. 
In this article, we prove the {\it arithmeticity} of the remaining 4 monodromy groups, for which, the parameters are $\alpha=(0,0,0,0)$; $\beta=(\frac{1}{3},\frac{1}{3},\frac{2}{3},\frac{2}{3}), (\frac{1}{3},\frac{2}{3},\frac{1}{6},\frac{5}{6}), (\frac{1}{4},\frac{1}{4},\frac{3}{4},\frac{3}{4}), (\frac{1}{3},\frac{2}{3},\frac{1}{4},\frac{3}{4})$; these are Examples 4, 8, 10, 11 of \cite{AvEvSZ, YYCE} and Examples \ref{arithmeticYYCE-4}, \ref{arithmeticYYCE-7}, \ref{arithmeticYYCE-5}, \ref{arithmeticYYCE-6} of Table \ref{table:calabiyau} (cf. \cite[Table 5.3]{SV}). In fact, we get the following theorem:
\begin{theorem}\label{maintheorem}
 The four hypergeometric monodromy groups associated to Calabi-Yau threefolds, for which, the parameters are $\alpha=(0,0,0,0)$; $\beta=(\frac{1}{3},\frac{1}{3},\frac{2}{3},\frac{2}{3})$, $(\frac{1}{3},\frac{2}{3},\frac{1}{6},\frac{5}{6}), (\frac{1}{4},\frac{1}{4},\frac{3}{4},\frac{3}{4}), (\frac{1}{3},\frac{2}{3},\frac{1}{4},\frac{3}{4})$, are arithmetic.
\end{theorem}
For a better reference, we list here the pairs of polynomials $f, g$ and the parameters $\alpha, \beta$, which correspond to the 14 hypergeometric monodromy groups associated to Calabi-Yau threefolds;  in the following list, \[\alpha=(0,0,0,0) \mbox{ i.e. } f(X)=(X-1)^4=X^4-4X^3+6X^2-4X+1:\]
{\renewcommand{\arraystretch}{1.3}   
{\tiny\begin{table}[h]
\addtolength{\tabcolsep}{0pt}
\caption{(cf. \cite[Table 5.3]{SV})}
\newcounter{rownum-3}
\setcounter{rownum-3}{0}
\centering
\begin{tabular}{ |c|  c|   c| c| c| c|}
\hline

  No. & $g(X)$ & $\beta$ & $f(X)-g(X)$ & Arithmetic\\ \hline
  
 { \refstepcounter{rownum-3}\arabic{rownum-3}\label{arithmeticYYCE-1}} &${ X^4-2X^3+3X^2-2X+1}$ & ${ 
\frac{1}{6}}$,${ \frac{1}{6}}$,${ \frac{5}{6}}$,${ \frac{5}{6} }$ & ${ -2X^3+3X^2-2X}$ & Yes, \cite{SV}\\ \hline

  { \refstepcounter{rownum-3}\arabic{rownum-3}*\label{brav1}} &${ X^4+4X^3+6X^2+4X+1}$ &${ \frac{1}{2}}$,${ \frac{1}{2}}$,${ \frac{1}{2}}$,${ \frac{1}{2}}$ &${ -8X^3-8X}$ &{ No, \cite{BT}} \\ \hline

 \refstepcounter{rownum-3}\arabic{rownum-3}\label{arithmeticYYCE-4} & $X^4+2X^3+3X^2+2X+1$& $\frac{1}{3}$,$\frac{1}{3}$,$\frac{2}{3}$,$\frac{2}{3}$ &$-6X^3+3X^2-6X$ & Yes \\ \hline
  
 {   \refstepcounter{rownum-3}\arabic{rownum-3}*\label{brav2}} &$ { X^4+3X^3+4X^2+3X+1}$ &$ { \frac{1}{2}}$,$ { \frac{1}{2}}$,$ { \frac{1}{3}}$,$ { \frac{2}{3}}$ &$ { -7X^3+2X^2-7X}$ & { No, \cite{BT}}	 \\ \hline
  
  \refstepcounter{rownum-3}\arabic{rownum-3}\label{arithmeticYYCE-5}  &$X^4+2X^2+1$  &$\frac{1}{4}$,$\frac{1}{4}$,$\frac{3}{4}$,$\frac{3}{4}$ &$-4X^3+4X^2-4X$ & Yes\\ \hline
  
 {   \refstepcounter{rownum-3}\arabic{rownum-3}*\label{brav3}} &$ { X^4+2X^3+2X^2+2X+1}$ &$ { \frac{1}{2}}$,$ { \frac{1}{2}}$,$ { \frac{1}{4}}$,$ { \frac{3}{4}}$ &$ { -6X^3+4X^2-6X}$& { No, \cite{BT}} \\ \hline
  
\refstepcounter{rownum-3}\arabic{rownum-3}\label{arithmeticYYCE-6}   &$X^4+X^3+2X^2+X+1$ &$\frac{1}{3}$,$\frac{2}{3}$,$\frac{1}{4}$,$\frac{3}{4}$ &$-5X^3+4X^2-5X$ &Yes \\ \hline
  
   {\refstepcounter{rownum-3}\arabic{rownum-3}*\label{brav4}}  & $ { X^4+X^3+X^2+X+1}$ &$ { 
\frac{1}{5}}$,$ { \frac{2}{5}}$,$ { \frac{3}{5}}$,$ { \frac{4}{5}}$ &$ { -5X^3+5X^2-5X}$ & { No, \cite{BT}}\\ \hline
  
 {\refstepcounter{rownum-3}\arabic{rownum-3}*\label{brav5}}  &$ { X^4+X^3+X+1}$ &$ { \frac{1}{2}}$,$ { \frac{1}{2}}$,$ { \frac{1}{6}}$,$ { \frac{5}{6}}$ &$ { -5X^3+6X^2-5X}$& { No, \cite{BT}}\\ \hline
  
 \refstepcounter{rownum-3}\arabic{rownum-3}\label{arithmeticYYCE-7} &$X^4+X^2+1$ &$\frac{1}{3}$,$\frac{2}{3}$,$\frac{1}{6}$,$\frac{5}{6}$ &$-4X^3+5X^2-4X$ &Yes\\ \hline
  
  {\refstepcounter{rownum-3}\arabic{rownum-3}\label{arithmeticYYCE-2}}  &${ X^4-X^3+2X^2-X+1}$ &${ \frac{1}{4}}$,${ \frac{3}{4}}$,${ \frac{1}{6}}$,${ \frac{5}{6}}$ &${ -3X^3+4X^2-3X}$ &{ Yes, \cite{SV}} \\ \hline
  
   { \refstepcounter{rownum-3}\arabic{rownum-3}*\label{brav6}} &$ { X^4+1}$ &$ { \frac{1}{8}}$,$ { \frac{3}{8}}$,$ { 
\frac{5}{8}}$,$ { \frac{7}{8}}$ &$ { -4X^3+6X^2-4X}$& { No, \cite{BT}}\\ \hline
  
  {\refstepcounter{rownum-3}\arabic{rownum-3}\label{arithmeticYYCE-3}} &${ X^4-X^3+X^2-X+1}$ &${ \frac{1}{10}}$,${ \frac{3}{10}}$,${ \frac{7}{10}}$,${ \frac{9}{10}}$ &${ -3X^3+5X^2-3X}$ &{ Yes, \cite{SV}}\\ \hline
  
   { \refstepcounter{rownum-3}\arabic{rownum-3}*\label{brav7}} &$ { X^4-X^2+1}$ &$ { \frac{1}{12}}$,$ { \frac{5}{12}}$,$ { \frac{7}{12}}$,$ { \frac{11}{12}}$ &$ { -4X^3+7X^2-4X}$ & { No, \cite{BT}}\\ \hline  
  \end{tabular}
\label{table:calabiyau}
\end{table}}

Therefore 7 of the 14 hypergeometric monodromy groups associated to Calabi-Yau threefolds, are arithmetic and other 7 are thin.  

For the pairs in Theorem \ref{maintheorem},  we explicitly compute, up to scalar multiples, the sypmplectic form $\Omega$ on $\Q^4$ and get a basis $\{\e_1,\e_2,\e_2^*,\e_1^*\}$ of $\Q^4$, satisfying the following: $\Omega(\e_i, \e_i*)\neq 0$ for $i=1,2$, and $\Omega(\e_1, \e_2)=\Omega(\e_1,\e_2^*)=\Omega(\e_2,\e_1^*)=\Omega(\e_2^*, \e_1^*)=0$. 

It can be shown easily that, with respect to the symplectic basis $\{\e_1,\e_2,\e_2^*,\e_1^*\}$ of $\Q^4$, the group of diagonal matrices in $\Sp_4(\Omega)$ form a maximal torus $\mathrm{T}$, the group of upper (resp. lower) triangular matrices in $\Sp_4(\Omega)$ form a Borel subgroup $\mathrm{B}$ (resp. $\mathrm{B}^-$, opposite to $\mathrm{B}$), and the group of unipotent upper (resp. lower) triangular matrices in $\Sp_4(\Omega)$ form the unipotent radical $\mathrm{U}$ (resp. $\mathrm{U}^-$, opposed to $\mathrm{U}$) of $\mathrm{B}$ (resp. $\mathrm{B}^-$). 

We now state a theorem of Tits (\cite{T}; cf. \cite[Theorem 7]{Ve2}) in a way we use it to prove Theorem \ref{maintheorem}.
\begin{theorem}\label{tits}
For every integer $n\geq 1$, the subgroup $\Gamma_n$ of $\Sp_4(\Z)$ generated by $n$-th power of the elements in $\mathrm{U}(\Z)$ and $\mathrm{U}^-(\Z)$ has finite index in $\Sp_4(\Z)$.
\end{theorem}
 We now have the following remarks:
\begin{rmk}\label{remark1} Note that if $\Gamma$ is a Zariski dense subgroup of $\Sp_4(\Z)$ such that $\Gamma\cap\mathrm{U}(\Z)$ is of finite index in $\mathrm{U}(\Z)$, then $\Gamma\cap\mathrm{U}^-(\Z)$ also has finite index in $\mathrm{U}^-(\Z)$ and there exists an integer $n\geq 1$ (big enough), so that the $n$-th power of elements in $\mathrm{U}(\Z)$ (resp. $\mathrm{U}^-(\Z)$) belong to $\Gamma\cap\mathrm{U}(\Z)$ (resp. $\Gamma\cap\mathrm{U}^-(\Z)$); and hence by Theorem \ref{tits}, $\Gamma$ has finite index in $\Sp_4(\Z)$.
\end{rmk}

\begin{rmk}\label{remark2} To prove Theorem \ref{maintheorem} using Remark \ref{remark1}, it is enough to show that $\Gamma\cap\mathrm{U}(\Z)$ is of finite index in $\mathrm{U}(\Z)$ (since the monodromy groups $\Gamma$ of Theorem \ref{maintheorem}  are Zariski dense subgroups of $\Sp_4(\Z)$).
\end{rmk}
Therefore our strategy is to show that the monodromy groups $\Gamma$ associated to the pairs in Theorem \ref{maintheorem}, intersects the group $\mathrm{U}(\Z)$ of unipotent upper (or lower) triangular matrices in $\Sp_4(\Omega)(\Z)$, in a finite index subgroup of $\mathrm{U}(\Z)$ i.e., $\Gamma\cap\mathrm{U}(\Z)$ is of finite index in $\mathrm{U}(\Z)$.
\begin{rmk}\label{venkataramana}
 It follows from \cite[Theorem 3.5]{Ve} that if $\Gamma$ is a Zariski dense subgroup of $\Sp_4(\Z)$ such that $\Gamma$ intersects the highest  and second highest root groups non-trivially, then  $\Gamma$ has finite index in $\Sp_4(\Z)$.
 
 Note that, once we show that the monodromy groups $\Gamma$ associated to the pairs in Theorem \ref{maintheorem}, intersects the group $\mathrm{U}(\Z)$ of unipotent upper (or lower) triangular matrices in $\Sp_4(\Omega)(\Z)$, in a finite index subgroup of $\mathrm{U}(\Z)$, it  follows automatically that $\Gamma$ also intersects the highest  and second highest root groups non-trivially; and hence the proof of Theorem \ref{maintheorem} also follows from \cite[Theorem 3.5]{Ve}.
\end{rmk}

\begin{rmk}
While the method of proof follows that of \cite{SV} closely, the actual computations in the present paper are much more complicated. In \cite{SV}, we use the following criterion to prove the arithmeticity for 3 (Examples \ref{arithmeticYYCE-1}, \ref{arithmeticYYCE-2}, \ref{arithmeticYYCE-3} of Table \ref{table:calabiyau}) of the 14 hypergeometric monodromy groups associated to Calabi-Yau threefolds: ``Let $A, B$ be the companion matrices of $f,g$ resp.,  $C=A^{-1}B$, and $v\in\Z^4$ such that $C(e_4)=e_4+v$, where $v$ is the linear combination of $e_1, e_2, e_3$; and $\{e_1, e_2, e_3, e_4\}$ is the standard basis of $\Q^4$ over $\Q$. If there exists an element $\gamma\in\Gamma=<A,B>$     such that the absolute value of $c$ ($\neq 0$) (the coefficient  of $e_4$ in $\gamma(v)$) is less or equal to $2$ and $\{v, \gamma^{-1}(v), \gamma(v)\}$ forms a linearly independent set, then the monodromy group $\Gamma$ is an arithmetic subgroup of $\Sp_4$.'' To prove the last statement, we show that the subgroup of $\Gamma$ 
generated 
by the reflections $C, \gamma^{-1}C\gamma, \gamma C\gamma^{-1}$ intersects (with respect to a suitable basis) 
the group $\mathrm{U}(\Z)$ of unipotent upper triangular matrices in $\Sp_4(\Z)$, in a finite index subgroup of $\mathrm{U}(\Z)$, by using the fact that the two matrices \[\begin{pmatrix}                                                                                                                                                                                                                                                                                                                                                                                                                                                                                                                                                                                                                                                                                                                                                                                                                                         1 &c\\                      
                                                                                                                                                                                                                                                                                                                                                                                                                                                                                                                                                                                                                                                                                                                                                                                                                0 &1                                                                                                                                                                                                                                            
                                                                                                                                                                                                                                                                                                                                                                                                                                                                                                                                                                                            \end{pmatrix}\mbox{ and } \begin{pmatrix}                                                                                                                                                                                                                                                                                                                                                                                                                           
                                                                                                                                                                                                                                                                                                                                                                                                              1 &0\\                                                                                                                                                                                                                                                                                                                                                                                                                                                                                                                                                                                                                                            
                                                                                                                                                                                          c &1                                                                                                                                                                                                                                                                                                                                                                                                                                                                                                                                                                                                                                                                                                                                                                                                                                        \end{pmatrix}\]for $1\leq\
\mid c\mid\leq 2$, generate a finite index subgroup of $\SL_2(\Z)$.

But for the 4 monodromy groups (Examples \ref{arithmeticYYCE-4}, \ref{arithmeticYYCE-5}, \ref{arithmeticYYCE-6}, \ref{arithmeticYYCE-7} of Table \ref{table:calabiyau}) mentioned in Theorem \ref{maintheorem}, it was not possible to get the element $\gamma\in\Gamma$ which satisfies the criterion of last para. Therefore we must make a good guess about the three reflections $P, Q, R$ (the conjugates of $C=A^{-1}B$) to be considered so that the technique of \cite{SV} can be applied. It is here that most of the difficulty lies. We make use of extensive computer calculations to confirm our guess about the three reflections to be chosen.
\end{rmk}

\section*{Acknowledgements}
I thank J\"org Hofmann for introducing me the program ``Maple'' to do the matrix computations; Maple was very helpful for the computations in this paper. I am grateful to Professor T. N. Venkataramana for introducing me the criterion to prove the arithmeticity. I thank him for his encouragement and constant support. I am also grateful to  Professor Duco van Straten and Professor Wadim Zudilin for their encouragement. I thank Institut f\"ur Mathematik, Johannes Gutenberg-Universit\"at for the postdoctoral fellowship, and for very pleasant hospitality. I also thank the referees for their valuable comments and suggestions.

\section{Proof of Theorem \ref{maintheorem}}
We will first compute the symplectic form $\Omega$ preserved by the monodromy group $\Gamma$, then show that there exists a basis $\{\e_1,\e_2,\e_2^*,\e_1^*\}$ of $\Q^4$, with respect to which, 
\[\Omega=\begin{pmatrix}
0 &0 &0 &\lambda_1\\
0 &0 &\lambda_2 &0\\
0 &-\lambda_2 &0 &0\\
-\lambda_1 &0 &0 &0
\end{pmatrix}\] where $\Omega(\e_i,\e_i^*)=\lambda_i\in\Q^{*},\ \forall 1\leq{i}\leq 2$.

It can be checked easily that the diagonal matrices in $\Sp_4(\Omega)$ form a maximal torus $\mathrm{T}$ i.e.,
\[\mathrm{T}=\left\{\begin{pmatrix}
t_1 &0 &0 &0\\
0 &t_2 &0 &0\\
0 &0 &t_2^{-1} &0\\
0 &0 &0 &t_1^{-1}
\end{pmatrix} : t_i\in\C^*,\quad \forall\ 1\leq i\leq 2\right\}\] is a maximal torus in $\Sp_4(\Omega)$. Once we fix a maximal torus $\mathrm{T}$ in $\Sp_4(\Omega)$, one may compute the root system $\Phi$ for $\Sp_4(\Omega)$. If we denote by $\mathbf{t}_i$, the character of $\mathrm{T}$ defined by
\[\begin{pmatrix}
t_1 &0 &0 &0\\
0 &t_2 &0 &0\\
0 &0 &t_2^{-1} &0\\
0 &0 &0 &t_1^{-1}
\end{pmatrix}\mapsto t_i,\qquad\mbox{for }i=1,2,\] then the roots are $$\Phi=\{\mathbf{t}_1^2, \mathbf{t}_1\mathbf{t}_2, \mathbf{t}_1\mathbf{t}_2^{-1}, \mathbf{t}_2^2, \mathbf{t}_1^{-2}, \mathbf{t}_1^{-1}\mathbf{t}_2^{-1}, \mathbf{t}_1^{-1}\mathbf{t}_2, \mathbf{t}_2^{-2}\}.$$ If we fix a set of simple roots $\Delta=\{\mathbf{t}_1\mathbf{t}_2^{-1}, \mathbf{t}_2^2\}$, then the set of positive roots $\Phi^+=\{\mathbf{t}_1^2, \mathbf{t}_1\mathbf{t}_2, \mathbf{t}_1\mathbf{t}_2^{-1}, \mathbf{t}_2^2\}$ and the group of upper (resp. lower) triangular matrices in $\Sp_4(\Omega)$ form a Borel subgroup $\mathrm{B}$ (resp. $\mathrm{B}^-$, opposite to $\mathrm{B}$), and the group of unipotent upper (resp. lower) triangular matrices in $\Sp_4(\Omega)$ form the unipotent radical $\mathrm{U}$ (resp. $\mathrm{U}^-$, opposed to $\mathrm{U}$) of $\mathrm{B}$ (resp. $\mathrm{B}^-$).

We now note the following remark:
\begin{rmk}\label{raghunathan}
 Since $\mathrm{U}$ is a nilpotent subgroup of $\GL_4(\R)$, it follows from \cite[Theorem 2.1]{Rag} that if $\Gamma\cap\mathrm{U}(\Z)$ is a Zariski dense subgroup of $\mathrm{U}$ then $\mathrm{U}/\Gamma\cap\mathrm{U}(\Z)$ is compact, and hence $\Gamma\cap\mathrm{U}(\Z)$ has finite index in $\mathrm{U}(\Z)$. 
\end{rmk}

Therefore, to show that the group $\Gamma$, with respect to the basis $\{\e_1,\e_2,\e_2^*,\e_1^*\}$ of $\Q^4$, intersects $\mathrm{U}(\Z)$ (resp. $\mathrm{U}^-(\Z)$) in a finite index subgroup of $\mathrm{U}(\Z)$ (resp. $\mathrm{U}^-(\Z)$), it is enough to show that $\Gamma$ contains non-trivial unipotent elements corresponding to each of the positive (resp. negative) roots (it shows that $\Gamma\cap\mathrm{U}(\Z)$ (resp. $\Gamma\cap\mathrm{U}^-(\Z)$) is Zariski dense in $\mathrm{U}$ (resp. $\mathrm{U}^-$)).

We will compute the symplectic form $\Omega$, the basis $\{\e_1,\e_2,\e_2^*,\e_1^*\}$ and the positive (or negative) root group elements $P, x, y, z$ (or $Q, x, y, z)\in\Gamma$ (for notation, see the proof below), and the proof  of arithmeticity of $\Gamma$ follows from Remark \ref{remark1} (cf. Remark \ref{venkataramana}). 

We will write the descriptions only for the pair $\alpha=(0,0,0,0)$, $\beta=(\frac{1}{3},\frac{2}{3},\frac{1}{4},\frac{3}{4})$ in Subsection \ref{calabi-yau11}, and write directly the form $\Omega$, the basis $\{\e_1,\e_2,\e_2^*,\e_1^*\}$, the elements $P, x, y, z$ (or $Q, x, y, z)\in\Gamma$ for other pairs, and the proof follows by the same descriptions as in Subsection \ref{calabi-yau11}.

\subsection{Arithmeticity of the monodromy group associated to the pair $\alpha=(0,0,0,0)$, $\beta=(\frac{1}{3},\frac{2}{3},\frac{1}{4},\frac{3}{4})$}\label{calabi-yau11}
This is Example 11 of \cite{AvEvSZ, YYCE} and Example \ref{arithmeticYYCE-6} of Table \ref{table:calabiyau} (cf. \cite[Table 5.3]{SV}). In this case $$f(X)=X^4-4X^3+6X^2-4X+1, \quad g(X)=X^4+X^3+2X^2+X+1;$$ and $f(X)-g(X)=-5X^3+4X^2-5X$.

Let $A$ and $B$  be the companion  matrices of $f(X)$ and  $g(X)$ resp., and let $C=A^{-1}B$. Then
{\scriptsize \[A=\begin{pmatrix}  \begin {array}{rrrr} 0&0&0&-1\\ \noalign{\medskip}1&0&0&4
\\ \noalign{\medskip}0&1&0&-6\\ \noalign{\medskip}0&0&1&4\end {array}
  \end{pmatrix},  B=\begin{pmatrix} \begin {array}{rrrr} 0&0&0&-1\\ \noalign{\medskip}1&0&0&-1
\\ \noalign{\medskip}0&1&0&-2\\ \noalign{\medskip}0&0&1&-1\end {array}
 \end{pmatrix}, C=A^{-1}B=\begin{pmatrix} \begin {array}{rrrr} 1&0&0&-5\\ \noalign{\medskip}0&1&0&4
\\ \noalign{\medskip}0&0&1&-5\\ \noalign{\medskip}0&0&0&1\end {array}
 \end{pmatrix}.\]}   Let $\Gamma=<A,B>$ be the  subgroup  of
$\SL_4(\Z)$ generated by $A$ and $B$. 
\subsection*{The invariant symplectic form} It follows form \cite{BH} that $\Gamma$ preserves a non-degenerate symplectic form $\Omega$ on $\Q^4$, which is integral on $\Z^4$ and the Zariski closure of $\Gamma$ is $\Sp_4(\Omega)$.  Let us denote $\Omega(v_1,v_2)$ by $v_1.v_2$, for any pairs of vectors $v_1,v_2\in\Q^4$. 

We compute here the form $\Omega$. Let $\{e_1,e_2,e_3,e_4\}$ be the standard basis of $\Q^4$ over $\Q$ and $v=-5e_1+4e_2-5e_3$ which is the last column vector of $C-\mathrm{I}_4$, where $\mathrm{I}_4$ is the $4\times 4$ identity matrix. Since $C$ preserves the form $\Omega$, for $1\leq i\leq 3$, we get 
\begin{align*}
e_i.e_4&=e_i.(-5e_1+4e_2-5e_3+e_4)\\
&=e_i.(v+e_4)\\
&=e_i.v+e_i.e_4.
\end{align*}
This implies that 
\begin{equation}\label{orthogonalv}
 e_i.v=0\qquad \mbox{for}\ 1\leq i\leq 3.
\end{equation}
That is, $v$ is $\Omega$- orthogonal to the vectors $e_1, e_2, e_3$ and $e_4.v\neq 0$ (since $\Omega$ is non-degenerate). Since $B$ preserves $\Omega$, we get
\begin{align*}
e_1.e_2=e_2.e_3=e_3.e_4&=e_4.(-e_1-e_2-2e_3-e_4)\\
&=e_1.e_4+e_2.e_4+2e_3.e_4.
\end{align*}
This implies that
\begin{equation}\label{equation2}
e_3.e_4=-e_1.e_4-e_2.e_4.
\end{equation}
It now follows from (\ref{orthogonalv}), (\ref{equation2}) and the invariance of $\Omega$ under $B$, that
\begin{equation}\label{equation3}
 e_1.e_3=\frac{4}{5}e_1.e_2=e_2.e_4
\end{equation} and
\begin{equation}\label{equation4}
 e_1.e_4=-\frac{9}{5}e_1.e_2.
\end{equation}

We now get from (\ref{equation2}), (\ref{equation3}) and (\ref{equation4}),  up to scalar multiples, the matrix form of {\scriptsize $$\Omega=\begin{pmatrix}
\begin {array}{rrrr} 0&1&4/5&-9/5\\ \noalign{\medskip}-1&0&1&4
/5\\ \noalign{\medskip}-4/5&-1&0&1\\ \noalign{\medskip}9/5&-4/5&-1&0
\end {array}
\end{pmatrix}.$$}
\subsection*{Proof of the arithmeticity of $\Gamma$} 
By an easy computation we get $\epsilon_1=e_1+e_3$, $\epsilon_2=-5e_1+4e_2-5e_3$, $\epsilon_2^*=-9e_1+4e_2-5e_3+4e_4$ and $\epsilon_1^*=e_1+e_2$ form a basis of $\Q^4$ over $\Q$, with respect to which 
{\scriptsize $$\Omega=\begin{pmatrix}\begin {array}{rrrr} 0&0&0&-4/5\\ \noalign{\medskip}0&0&{
\frac {144}{5}}&0\\ \noalign{\medskip}0&-{\frac {144}{5}}&0&0
\\ \noalign{\medskip}4/5&0&0&0\end {array}
\end{pmatrix}.$$}
Let $$P=C=A^{-1}B,\quad  Q=B^{-3}CB^3,\quad R=B^3CB^{-3}.$$   It can be
checked     easily    that     with    respect     to     the    basis
$\{\epsilon_1,\epsilon_2,\epsilon_2^*,\epsilon_1^*\}$,  the $P, Q, R$ have, resp., the matrix form
{\scriptsize \[\begin{pmatrix}\begin {array}{rrrr} 1&0&0&0\\ \noalign{\medskip}0&1&4&0
\\ \noalign{\medskip}0&0&1&0\\ \noalign{\medskip}0&0&0&1\end {array}
 \end{pmatrix}, \qquad  \begin{pmatrix} \begin {array}{rrrr} 1&0&0&0\\ \noalign{\medskip}0&1&0&0
\\ \noalign{\medskip}0&-4&1&0\\ \noalign{\medskip}0&0&0&1\end {array}
   \end{pmatrix}, \qquad \begin{pmatrix}  \begin {array}{rrrr} 1&-72&0&-36\\ \noalign{\medskip}0&1&0&0
\\ \noalign{\medskip}0&-4&1&-2\\ \noalign{\medskip}0&0&0&1\end {array}
 \end{pmatrix}.\]}  A  computation  shows  that  if
\[E=Q^{-1}R,\ F=Q^{-1}RQP^{-1},\ G=[E,F]=EFE^{-1}F^{-1},\ x=[G,E],\] \[y=x^{-36}E^{1152},\ u=G^{82944}y^{1152},\ z=u^{-1152}x^{63403237638144},\] then

{\scriptsize\[E=\begin{pmatrix}
\begin {array}{rrrr} 1&-72&0&-36\\ \noalign{\medskip}0&1&0&0
\\ \noalign{\medskip}0&0&1&-2\\ \noalign{\medskip}0&0&0&1\end {array}
\end{pmatrix},\quad F=\begin{pmatrix}
\begin {array}{rrrr} 1&-72&288&-36\\ \noalign{\medskip}0&1&-4&0
\\ \noalign{\medskip}0&-4&17&-2\\ \noalign{\medskip}0&0&0&1
\end {array}
\end{pmatrix},\] \[G=\begin{pmatrix}\begin {array}{rrrr} 1&1152&288&1728\\ \noalign{\medskip}0&1&0
&-8\\ \noalign{\medskip}0&0&1&32\\ \noalign{\medskip}0&0&0&1\end {array}
\end{pmatrix},\quad u=\begin{pmatrix}
             \begin {array}{rrrr} 1&0&23887872&-63403237638144
\\ \noalign{\medskip}0&1&0&-663552\\ \noalign{\medskip}0&0&1&0
\\ \noalign{\medskip}0&0&0&1\end {array}
            \end{pmatrix},\]} {\scriptsize\[x=\begin{pmatrix} \begin {array}{rrrr} 1&0&0&-1152\\ \noalign{\medskip}0&1&0&0
\\ \noalign{\medskip}0&0&1&0\\ \noalign{\medskip}0&0&0&1\end {array}
\end{pmatrix}, \quad y=\begin{pmatrix} \begin {array}{rrrr} 1&-82944&0&0\\ \noalign{\medskip}0&1&0&0
\\ \noalign{\medskip}0&0&1&-2304\\ \noalign{\medskip}0&0&0&1
\end {array}
\end{pmatrix},\]} {\scriptsize\[\quad z=\begin{pmatrix}\begin{array}{rrrr} 1&0&-27518828544&0\\ \noalign{\medskip}0&
1&0&764411904\\ \noalign{\medskip}0&0&1&0\\ \noalign{\medskip}0&0&0&1
\end {array}\end{pmatrix}.\]}  The  elements $P, x, y, z$ patently generate  the positive  root  groups  of
$\Sp_4(\Omega)$. Therefore if $\mathrm{B}$ is the Borel subgroup of
$\Sp_4(\Omega)(\Q)$ preserving the full flag
\[\{0\}\subset\Q\e_1\subset\Q\epsilon_1\oplus\Q\epsilon_2
\subset\Q\epsilon_1
\oplus\Q\epsilon_2\oplus\Q\epsilon_2^*\subset\Q\epsilon_1
\oplus\Q\epsilon_2\oplus\Q\epsilon_2^*\oplus\Q\epsilon_1^*=\Q^4\]   and
$\mathrm{U}$ its unipotent radical, then $\Gamma\cap\mathrm{U}(\Z)$ is
of finite index in $\mathrm{U}(\Z)$. Since $\Gamma$ is Zariski dense in
$\Sp_4(\Omega)$ by \cite{BH}, it follows  from Remark \ref{remark1} (cf. Remark \ref{venkataramana})
that $\Gamma$ is an arithmetic subgroup of $\Sp_4(\Omega)(\Z)$.\qed

\subsection{Arithmeticity of the monodromy group associated to the pair $\alpha=(0,0,0,0)$, $\beta=(\frac{1}{3},\frac{1}{3},\frac{2}{3},\frac{2}{3})$}\label{calabi-yau4}
This is Example 4 of \cite{AvEvSZ, YYCE} and Example \ref{arithmeticYYCE-4} of Table \ref{table:calabiyau} (cf. \cite[Table 5.3]{SV}). In this case $$f(X)=X^4-4X^3+6X^2-4X+1, \quad g(X)=X^4+2X^3+3X^2+2X+1;$$ and $f(X)-g(X)=-6X^3+3X^2-6X$.

Let $A$ and $B$  be the companion  matrices of $f(X)$ and  $g(X)$ resp., and let $C=A^{-1}B$. Then
{\scriptsize \[A=\begin{pmatrix}\begin {array}{rrrr} 0&0&0&-1\\ \noalign{\medskip}1&0&0&4
\\ \noalign{\medskip}0&1&0&-6\\ \noalign{\medskip}0&0&1&4\end {array}
  \end{pmatrix},  B=\begin{pmatrix}\begin {array}{rrrr} 0&0&0&-1\\ \noalign{\medskip}1&0&0&-2
\\ \noalign{\medskip}0&1&0&-3\\ \noalign{\medskip}0&0&1&-2\end {array}
 \end{pmatrix}, C=A^{-1}B=\begin{pmatrix}\begin {array}{rrrr} 1&0&0&-6\\ \noalign{\medskip}0&1&0&3
\\ \noalign{\medskip}0&0&1&-6\\ \noalign{\medskip}0&0&0&1\end {array}
 \end{pmatrix}.\]}   Let $\Gamma=<A,B>$ be the  subgroup  of
$\SL_4(\Z)$ generated by $A$ and $B$. 

\subsection*{The invariant symplectic form} Using the same method as in Subsection \ref{calabi-yau11}, we get the matrix form of
{\scriptsize\[\Omega=\begin{pmatrix}
         \begin {array}{rrrr} 0&1&1/2&-3\\ \noalign{\medskip}-1&0&1&1/2
\\ \noalign{\medskip}-1/2&-1&0&1\\ \noalign{\medskip}3&-1/2&-1&0
\end {array}
\end{pmatrix}.\]}
\subsection*{Proof of the arithmeticity of $\Gamma$}
By an easy computation we get $\epsilon_1=e_1+e_3$, $\epsilon_2=-6e_1+3e_2-6e_3$, $\epsilon_2^*=24e_1+12e_2+12e_3-3e_4$ and $\epsilon_1^*=e_1+2e_2$ form a basis of $\Q^4$ over $\Q$, with respect to which 
{\scriptsize $$\Omega=\begin{pmatrix}\begin {array}{rrrr} 0&0&0&-1/2\\ \noalign{\medskip}0&0&-{
\frac {81}{2}}&0\\ \noalign{\medskip}0&{\frac {81}{2}}&0&0
\\ \noalign{\medskip}1/2&0&0&0\end {array}
\end{pmatrix}.$$}
Let $$P=C=A^{-1}B,\quad  Q=B^{-4}CB^4,\quad R=A^{-1}CA.$$   It can be
checked     easily    that     with    respect     to     the    basis
$\{\epsilon_1,\epsilon_2,\epsilon_2^*,\epsilon_1^*\}$,  the $P, Q, R$ have, resp., the matrix form
{\scriptsize \[\begin{pmatrix} \begin {array}{rrrr} 1&0&0&0\\ \noalign{\medskip}0&1&-3&0
\\ \noalign{\medskip}0&0&1&0\\ \noalign{\medskip}0&0&0&1\end {array}
 \end{pmatrix}, \qquad  \begin{pmatrix} \begin {array}{rrrr} 1&0&0&0\\ \noalign{\medskip}0&1&0&0
\\ \noalign{\medskip}0&3&1&0\\ \noalign{\medskip}0&0&0&1\end {array}
   \end{pmatrix}, \qquad \begin{pmatrix}\begin {array}{rrrr} 1&0&0&0\\ \noalign{\medskip}0&1&0&0
\\ \noalign{\medskip}-2&12&1&0\\ \noalign{\medskip}27&-162&0&1
\end {array}
 \end{pmatrix}.\]}  A  computation  shows  that  if
\[G=Q^{-4}R, \ H=PGP^{-1}, \ x=[G,H],\ y=x^{27}G^{1944},\ E=H^{1944}x^{27},\] \[\ F=E^{-1}y, \ z=F^{-1944}x^{3673320192},\] then
{\scriptsize\[G=\begin{pmatrix}\begin {array}{rrrr} 1&0&0&0\\ \noalign{\medskip}0&1&0&0
\\ \noalign{\medskip}-2&0&1&0\\ \noalign{\medskip}27&-162&0&1
\end {array}
\end{pmatrix},\qquad H=\begin{pmatrix}\begin {array}{rrrr} 1&0&0&0\\ \noalign{\medskip}6&1&0&0
\\ \noalign{\medskip}-2&0&1&0\\ \noalign{\medskip}27&-162&-486&1
\end {array}
\end{pmatrix},\]} {\scriptsize\[x=\begin{pmatrix} \begin {array}{rrrr} 1&0&0&0\\ \noalign{\medskip}0&1&0&0
\\ \noalign{\medskip}0&0&1&0\\ \noalign{\medskip}-1944&0&0&1
\end {array}
\end{pmatrix}, \quad y=\begin{pmatrix} \begin {array}{rrrr} 1&0&0&0\\ \noalign{\medskip}0&1&0&0
\\ \noalign{\medskip}-3888&0&1&0\\ \noalign{\medskip}0&-314928&0&1
\end {array}
\end{pmatrix},\] \[ E=\begin{pmatrix}\begin {array}{rrrr} 1&0&0&0\\ \noalign{\medskip}11664&1&0&0
\\ \noalign{\medskip}-3888&0&1&0\\ \noalign{\medskip}0&-314928&-944784
&1\end {array} \end{pmatrix},\ F=\begin{pmatrix}\begin {array}{rrrr} 1&0&0&0\\ \noalign{\medskip}-11664&1&0&0
\\ \noalign{\medskip}0&0&1&0\\ \noalign{\medskip}-3673320192&0&944784&
1\end {array}\end{pmatrix},\] \[z=\begin{pmatrix} \begin {array}{rrrr} 1&0&0&0\\ \noalign{\medskip}22674816&1&0&0
\\ \noalign{\medskip}0&0&1&0\\ \noalign{\medskip}0&0&-1836660096&1
\end {array}
\end{pmatrix}.\]}  The  elements $Q, x, y, z$ patently generate  the negative  root  groups  of
$\Sp_4(\Omega)$, therefore $\Gamma$ intersects the group $\mathrm{U}^-(\Z)$ of unipotent lower triangular matrices in $\Sp_4(\Omega)(\Z)$, in a finite index subgroup of $\mathrm{U}^-(\Z)$ i.e., $\Gamma\cap\mathrm{U}^-(\Z)$ is of finite index in $\mathrm{U}^-(\Z)$. Since $\Gamma$ is Zariski dense in
$\Sp_4(\Omega)$ by \cite{BH}, it follows from Remark \ref{remark1} (cf. Remark \ref{venkataramana}) that $\Gamma$ is an arithmetic subgroup of $\Sp_4(\Omega)(\Z)$.\qed

\subsection{Arithmeticity of the monodromy group associated to the pair $\alpha=(0,0,0,0)$, $\beta=(\frac{1}{3},\frac{2}{3},\frac{1}{6},\frac{5}{6})$}\label{calabi-yau8}
This is Example 8 of \cite{AvEvSZ, YYCE} and Example \ref{arithmeticYYCE-7} of Table \ref{table:calabiyau} (cf. \cite[Table 5.3]{SV}). In this case $$f(X)=X^4-4X^3+6X^2-4X+1, \quad g(X)=X^4+X^2+1;$$ and $f(X)-g(X)=-4X^3+5X^2-4X$.

Let $A$ and $B$  be the companion  matrices of $f(X)$ and  $g(X)$ resp., and let $C=A^{-1}B$. Then
{\scriptsize \[A=\begin{pmatrix}\begin {array}{rrrr} 0&0&0&-1\\ \noalign{\medskip}1&0&0&4
\\ \noalign{\medskip}0&1&0&-6\\ \noalign{\medskip}0&0&1&4\end {array}
  \end{pmatrix},  B=\begin{pmatrix}\begin {array}{rrrr} 0&0&0&-1\\ \noalign{\medskip}1&0&0&0
\\ \noalign{\medskip}0&1&0&-1\\ \noalign{\medskip}0&0&1&0\end {array}
 \end{pmatrix}, C=A^{-1}B=\begin{pmatrix}\begin {array}{rrrr} 1&0&0&-4\\ \noalign{\medskip}0&1&0&5
\\ \noalign{\medskip}0&0&1&-4\\ \noalign{\medskip}0&0&0&1\end {array}
 \end{pmatrix}.\]}   Let $\Gamma=<A,B>$ be the  subgroup  of
$\SL_4(\Z)$ generated by $A$ and $B$. 

\subsection*{The invariant symplectic form} Using the same method as in Subsection \ref{calabi-yau11}, we get the matrix form of
{\scriptsize\[\Omega=\begin{pmatrix}
        \begin {array}{rrrr} 0&1&5/4&0\\ \noalign{\medskip}-1&0&1&5/4
\\ \noalign{\medskip}-5/4&-1&0&1\\ \noalign{\medskip}0&-5/4&-1&0
\end {array}
\end{pmatrix}.\]}
\subsection*{Proof of the arithmeticity of $\Gamma$}
By an easy computation we get $\epsilon_1=e_2$, $\epsilon_2=-4e_1+5e_2-4e_3$, $\epsilon_2^*=5e_1+4e_4$ and $\epsilon_1^*=e_1$ form a basis of $\Q^4$ over $\Q$, with respect to which 
{\scriptsize $$\Omega=\begin{pmatrix}\begin {array}{rrrr} 0&0&0&-1\\ \noalign{\medskip}0&0&9&0
\\ \noalign{\medskip}0&-9&0&0\\ \noalign{\medskip}1&0&0&0\end {array}
\end{pmatrix}.$$}
Let $$P=C=A^{-1}B,\quad  Q=B^{-1}CB,\quad R=B^3CB^{-3},\quad S=P^{25}R^{-4}.$$   It can be
checked     easily    that     with    respect     to     the    basis
$\{\epsilon_1,\epsilon_2,\epsilon_2^*,\epsilon_1^*\}$,  the $P, Q, R, S$ have, resp., the matrix form
{\scriptsize \[\begin{pmatrix} \begin {array}{rrrr} 1&0&0&0\\ \noalign{\medskip}0&1&4&0
\\ \noalign{\medskip}0&0&1&0\\ \noalign{\medskip}0&0&0&1\end {array}
 \end{pmatrix},\begin{pmatrix} \begin {array}{rrrr} 1&0&0&0\\ \noalign{\medskip}0&1&0&0
\\ \noalign{\medskip}0&-4&1&0\\ \noalign{\medskip}0&0&0&1\end {array}
   \end{pmatrix}, \begin{pmatrix}\begin {array}{rrrr} 1&0&-45&-9\\ \noalign{\medskip}0&1&25&5
\\ \noalign{\medskip}0&0&1&0\\ \noalign{\medskip}0&0&0&1\end {array}
 \end{pmatrix}, \begin{pmatrix}\begin {array}{rrrr} 1&0&180&36\\ \noalign{\medskip}0&1&0&-20
\\ \noalign{\medskip}0&0&1&0\\ \noalign{\medskip}0&0&0&1\end {array}
\end{pmatrix}.\]}  A  computation  shows  that  if
\[G=QSQ^{-1}, \ x=[S,G],\ y=S^{28800}x^{-36},\ H=G^{5184000}y^{-180},\] \[z=H^{28800}x^{-386983526586624000},\] then
{\scriptsize\[G=\begin{pmatrix}\begin {array}{rrrr} 1&720&180&36\\ \noalign{\medskip}0&1&0&-
20\\ \noalign{\medskip}0&0&1&80\\ \noalign{\medskip}0&0&0&1
\end {array}
\end{pmatrix},\qquad x=\begin{pmatrix}\begin {array}{rrrr} 1&0&0&28800\\ \noalign{\medskip}0&1&0&0
\\ \noalign{\medskip}0&0&1&0\\ \noalign{\medskip}0&0&0&1\end {array}
\end{pmatrix},\]} {\scriptsize\[y=\begin{pmatrix} \begin {array}{rrrr} 1&0&5184000&0\\ \noalign{\medskip}0&1&0&-
576000\\ \noalign{\medskip}0&0&1&0\\ \noalign{\medskip}0&0&0&1
\end {array}
\end{pmatrix}, H=\begin{pmatrix}\begin {array}{rrrr} 1&3732480000&0&386983526586624000
\\ \noalign{\medskip}0&1&0&0\\ \noalign{\medskip}0&0&1&414720000
\\ \noalign{\medskip}0&0&0&1\end {array}
\end{pmatrix},\] \[z=\begin{pmatrix}\begin {array}{rrrr} 1&107495424000000&0&0
\\ \noalign{\medskip}0&1&0&0\\ \noalign{\medskip}0&0&1&11943936000000
\\ \noalign{\medskip}0&0&0&1\end {array}
\end{pmatrix}.\]}  The  elements $P, x, y, z$ patently generate  the positive  root  groups  of
$\Sp_4(\Omega)$. Since $\Gamma$ is Zariski dense in
$\Sp_4(\Omega)$ by \cite{BH}, it follows  from Remark \ref{remark1} (cf. Remark \ref{venkataramana})
that $\Gamma$ is an arithmetic subgroup of $\Sp_4(\Omega)(\Z)$.\qed

\subsection{Arithmeticity of the monodromy group associated to the pair $\alpha=(0,0,0,0)$, $\beta=(\frac{1}{4},\frac{1}{4},\frac{3}{4},\frac{3}{4})$}\label{calabi-yau10}
This is Example 10 of \cite{AvEvSZ, YYCE} and Example \ref{arithmeticYYCE-5} of Table \ref{table:calabiyau} (cf. \cite[Table 5.3]{SV}). In this case $$f(X)=X^4-4X^3+6X^2-4X+1, \quad g(X)=X^4+2X^2+1;$$ and $f(X)-g(X)=-4X^3+4X^2-4X$.

Let $A$ and $B$  be the companion  matrices of $f(X)$ and  $g(X)$ resp., and let $C=A^{-1}B$. Then
{\scriptsize \[A=\begin{pmatrix}\begin {array}{rrrr} 0&0&0&-1\\ \noalign{\medskip}1&0&0&4
\\ \noalign{\medskip}0&1&0&-6\\ \noalign{\medskip}0&0&1&4\end {array}
  \end{pmatrix},  B=\begin{pmatrix}\begin {array}{rrrr} 0&0&0&-1\\ \noalign{\medskip}1&0&0&0
\\ \noalign{\medskip}0&1&0&-2\\ \noalign{\medskip}0&0&1&0\end {array}
 \end{pmatrix}, C=A^{-1}B=\begin{pmatrix}\begin {array}{rrrr} 1&0&0&-4\\ \noalign{\medskip}0&1&0&4
\\ \noalign{\medskip}0&0&1&-4\\ \noalign{\medskip}0&0&0&1\end {array}
 \end{pmatrix}.\]}   Let $\Gamma=<A,B>$ be the  subgroup  of
$\SL_4(\Z)$ generated by $A$ and $B$. 

\subsection*{The invariant symplectic form} Using the same method as in Subsection \ref{calabi-yau11}, we get the matrix form of
{\scriptsize\[\Omega=\begin{pmatrix}
        \begin {array}{rrrr} 0&1&1&-1\\ \noalign{\medskip}-1&0&1&1
\\ \noalign{\medskip}-1&-1&0&1\\ \noalign{\medskip}1&-1&-1&0
\end {array}
\end{pmatrix}.\]}
\subsection*{Proof of the arithmeticity of $\Gamma$}
By an easy computation we get $\epsilon_1=e_2$, $\epsilon_2=-e_1+e_2-e_3$, $\epsilon_2^*=3e_1+e_2+2e_3+e_4$ and $\epsilon_1^*=3e_1+2e_3$ form a basis of $\Q^4$ over $\Q$, with respect to which 
{\scriptsize $$\Omega=\begin{pmatrix}\begin {array}{rrrr} 0&0&0&-1\\ \noalign{\medskip}0&0&1&0
\\ \noalign{\medskip}0&-1&0&0\\ \noalign{\medskip}1&0&0&0\end {array}
\end{pmatrix}.$$}
Let $$P=C=A^{-1}B,\quad  Q=B^{-5}CB^5,\quad R=B^3CB^{-3},\quad S=P^{-1}RPQ^{-1}.$$   It can be
checked     easily    that     with    respect     to     the    basis
$\{\epsilon_1,\epsilon_2,\epsilon_2^*,\epsilon_1^*\}$,  the $P, Q, R, S$ have, resp., the matrix form
{\scriptsize \[\begin{pmatrix} \begin {array}{rrrr} 1&0&0&0\\ \noalign{\medskip}0&1&4&0
\\ \noalign{\medskip}0&0&1&0\\ \noalign{\medskip}0&0&0&1\end {array}
 \end{pmatrix},\begin{pmatrix} \begin {array}{rrrr} 1&0&0&0\\ \noalign{\medskip}0&1&0&0
\\ \noalign{\medskip}0&-4&1&0\\ \noalign{\medskip}0&0&0&1\end {array}
   \end{pmatrix}, \begin{pmatrix}\begin {array}{rrrr} 1&16&-64&-64\\ \noalign{\medskip}0&-15&64
&64\\ \noalign{\medskip}0&-4&17&16\\ \noalign{\medskip}0&0&0&1
\end {array}
 \end{pmatrix}, \begin{pmatrix}\begin {array}{rrrr} 1&16&0&-64\\ \noalign{\medskip}0&1&0&0
\\ \noalign{\medskip}0&0&1&16\\ \noalign{\medskip}0&0&0&1\end {array}
\end{pmatrix}.\]}  A  computation  shows  that  if
\[E=PSP^{-1},\ x=[S,E],\ y=S^{2048}x^{64},\ G=E^{2048}x^{64},\ H=G^{-1}y,\] \[z=H^{-2048}x^{4294967296},\] then
{\scriptsize\[E=\begin{pmatrix}\begin {array}{rrrr} 1&16&-64&-64\\ \noalign{\medskip}0&1&0&64
\\ \noalign{\medskip}0&0&1&16\\ \noalign{\medskip}0&0&0&1\end {array}
\end{pmatrix},\quad x=\begin{pmatrix}\begin {array}{rrrr} 1&0&0&2048\\ \noalign{\medskip}0&1&0&0
\\ \noalign{\medskip}0&0&1&0\\ \noalign{\medskip}0&0&0&1\end {array}
\end{pmatrix},\]} {\scriptsize\[y=\begin{pmatrix} \begin {array}{rrrr} 1&32768&0&0\\ \noalign{\medskip}0&1&0&0
\\ \noalign{\medskip}0&0&1&32768\\ \noalign{\medskip}0&0&0&1
\end {array}
\end{pmatrix},\quad G=\begin{pmatrix}\begin {array}{rrrr} 1&32768&-131072&0\\ \noalign{\medskip}0&1
&0&131072\\ \noalign{\medskip}0&0&1&32768\\ \noalign{\medskip}0&0&0&1
\end {array}
\end{pmatrix},\] \[H=\begin{pmatrix}\begin {array}{rrrr} 1&0&131072&4294967296
\\ \noalign{\medskip}0&1&0&-131072\\ \noalign{\medskip}0&0&1&0
\\ \noalign{\medskip}0&0&0&1\end {array}
\end{pmatrix},\quad z=\begin{pmatrix}\begin {array}{rrrr} 1&0&-268435456&0\\ \noalign{\medskip}0&1&0
&268435456\\ \noalign{\medskip}0&0&1&0\\ \noalign{\medskip}0&0&0&1
\end {array}
\end{pmatrix}.\]}  The  elements $P, x, y, z$ patently generate  the positive  root  groups  of
$\Sp_4(\Omega)$. Since $\Gamma$ is Zariski dense in
$\Sp_4(\Omega)$ by \cite{BH}, it follows  from Remark \ref{remark1} (cf. Remark \ref{venkataramana})
that $\Gamma$ is an arithmetic subgroup of $\Sp_4(\Omega)(\Z)$.\qed

\end{document}